\documentclass[aap,MSNbibl,seceqn,dvips]{arximspdf}

\doi{10.1214/13-AAP949} 
\volume{24}
\issue{4}
\pubyear{2014}
\firstpage{1375}
\lastpage{1395}

\makeatletter
\newcommand{\eqref}[1]{(\ref{#1})}
\newcommand{\Q}{\mathbb{Q}}
\newcommand{\R}{\mathbb{R}}
\newcommand{\bmm}[1]{\mathbf{#1}}

\renewcommand{\epsilon}{\varepsilon}

\newtheorem{theorem}{Theorem}[section]
\newtheorem{lemma}[theorem]{Lemma}
\newtheorem{proposition}[theorem]{Proposition}
\newtheorem{corollary}[theorem]{Corollary}

\newproclaim{rem}[theorem]{Remark}
\newproclaim{remark}[theorem]{Remark}
\newproclaim{model}{Model}
\newproclaim{definition}[theorem]{Definition}

\makeatother

\begin{document}
\begin{frontmatter}

\title{On the existence of accessible paths in various models of
fitness landscapes}
\runtitle{Accessible paths in fitness landscapes}

\begin{aug}
\author{\fnms{Peter} \snm{Hegarty}\corref{}\ead[label=e1]{hegarty@chalmers.se}}
\and
\author{\fnms{Anders} \snm{Martinsson}\ead[label=e2]{andemar@chalmers.se}}
\runauthor{P. Hegarty and A. Martinsson}
\affiliation{Chalmers University of Technology and University of Gothenburg}
\address{Department of Mathematical Sciences\\
Chalmers University Of Technology and\\
\quad University of Gothenburg\\
41296 Gothenburg\\
Sweden\\
\printead{e1}\\
\phantom{E-mail:\ }\printead*{e2}}
\end{aug}

\received{\smonth{10} \syear{2012}}
\revised{\smonth{7} \syear{2013}}

%
\begin{abstract} We present rigorous mathematical analyses of a number
of well-known mathematical models for genetic mutations. In these
models, the genome is represented by a vertex of the $n$-dimensional
binary hypercube, for some~$n$, a mutation involves the flipping of
a single bit, and each vertex is assigned a real number, called its
fitness, according to some rules. Our main concern is with the issue of
existence of (selectively) accessible paths; that is, monotonic paths
in the hypercube along which fitness is always increasing. Our main
results resolve open questions about three such models, which in the
biophysics literature are known as house of cards (HoC), constrained
house of cards (CHoC) and rough Mount Fuji (RMF). We prove that the
probability of there being at least one accessible path from the
all-zeroes node $\bmm{v}^0$ to the all-ones node $\bmm{v}^1$ tends
respectively to 0, 1 and 1, as $n$ tends to infinity. A crucial idea is
the introduction of a generalization of the CHoC model, in which the
fitness of $\bmm{v}^0$ is set to some $\alpha= \alpha_n \in[0, 1]$.
We prove that there is a very sharp threshold at $\alpha_n = \frac
{\ln
n}{n}$ for the existence of accessible paths from $\bmm{v}^0$ to~$\bmm{v}^1$.
As a corollary we prove significant concentration, for $\alpha$ below
the threshold, of the number of accessible paths about the expected
value (the precise statement is technical; see Corollary \ref{cor:aHoCconc}).
In the case of RMF, we prove that the probability of accessible paths
from $\bmm{v}^0$ to $\bmm{v}^1$
existing tends to $1$ provided the drift parameter $\theta= \theta_n$
satisfies $n\theta_n \rightarrow\infty$, and for any fitness
distribution which is continuous on its support and whose support is connected.
\end{abstract}




\end{frontmatter}
\setcounter{section}{-1}
\section{Notation}
Throughout this paper, $\mathbb{Q}_n$ will denote the \emph{directed}
$n$-dimensional binary hypercube. This is the directed graph whose nodes
are all binary strings of length $n$, with an edge between any
pair of nodes that differ in exactly one bit, the edge being
always directed toward the node with the greater number of ones.

Let $g, h \dvtx \mathbb{N} \rightarrow\mathbb{R}_{+}$ be any two
functions. We will employ the following notation throughout, all of
which is quite standard:
\begin{longlist}[(vii)]
\item[(i)] $g(n) \sim h(n)$ means that $\lim_{n \rightarrow\infty}
\frac
{g(n)}{h(n)} = 1$;

\item[(ii)] $g(n) \lesssim h(n)$ means that $\limsup_{n \rightarrow
\infty} \frac{g(n)}{h(n)} \leq1$;

\item[(iii)] $g(n) \gtrsim h(n)$ means that $h(n) \lesssim g(n)$;

\item[(iv)] $g(n) = O(h(n))$ means that $\limsup_{n \rightarrow\infty}
\frac{g(n)}{h(n)} < \infty$;

\item[(v)] $g(n) = \Omega(h(n))$ means that $h(n) = O(g(n))$;

\item[(vi)] $g(n) = \Theta(h(n))$ means that both $g(n) = O(h(n))$ and
$h(n) = O(g(n))$ hold;

\item[(vii)] $g(n) = o(h(n))$ means that $\lim_{n \rightarrow\infty}
\frac{g(n)}{h(n)} = 0$.
\end{longlist}

Now suppose instead that $(g(n))_{n=1}^{\infty}, (h(n))_{n=1}^{\infty}$
are two
sequences of random variables. We write $g(n) \sim h(n)$ if,
for all $\varepsilon_1, \varepsilon_2 > 0$ and $n$ sufficiently large,
%
\begin{equation}
\mathbb{P} \biggl( 1-\varepsilon_1 < \frac{g(n)}{h(n)} < 1 +
\varepsilon _1 \biggr) > 1-\varepsilon_2.
\end{equation}
Similarly, we write $g(n) \gtrsim h(n)$ if, for all $\varepsilon_1,
\varepsilon_2 > 0$ and $n$
sufficiently large,
%
\begin{equation}
\mathbb{P} \biggl( \frac{g(n)}{h(n)} > 1 - \varepsilon_1 \biggr) > 1
- \varepsilon_2.
\end{equation}

\section{Introduction}
In many basic mathematical models of genetic mutations, the genome is
represented as a node of the directed $n$-dimensional binary hypercube
$\mathbb{Q}_{n}$, and each mutation involves the flipping of a single
bit from
$0$ (the ``wild'' state) to $1$ (the ``mutant'' state), hence displacement
along an edge of $\mathbb{Q}_n$. Each node $v \in\mathbb{Q}_n$ is
assigned a
real number $f(v)$, called its \emph{fitness}. The fitness of a node is
not a
constant, but is drawn from some probability distribution specified by the
model. This distribution may vary from node to node in more or less complicated
ways, depending on the model. Basically, however, evolution is
considered as
favoring mutational pathways which, on average, lead to higher fitness. A
fundamental concept in this regard is the following (see \cite{W2,W1,FKVK}):

\begin{definition}\label{de1.1}
Let $f \dvtx \mathbb{Q}_n \rightarrow\mathbb{R}$ be
a fitness function. A \emph{(selectively) accessible path} in $\mathbb
{Q}_n$ is a path
%
\begin{equation}
v_0 \rightarrow v_1 \rightarrow\cdots\rightarrow
v_{k-1} \rightarrow v_k,
\end{equation}
such that $f(v_i) > f(v_{i-1})$ for $i = 1, \ldots, k$.
\end{definition}

Let $\bmm{v}^{0} = (0, 0, \ldots, 0)$, $\bmm{v}^{1} = (1, 1, \ldots, 1)$
denote the all-zeroes and all-ones vertices in $\mathbb{Q}_n$. A basic
question in such models is whether accessible paths from $\bmm{v}^0$ to~$\bmm{v}^1$ exist or not with high probability. For the remainder of
this paper, unless explicitly stated otherwise, the words ``accessible
path'' will always refer to such a path which starts at $\bmm{v}^0$ and
ends at $\bmm{v}^1$. In fact, it will only be in the proof of
Proposition~\ref{prop:bootstrap} that we will need to consider
accessible paths with other start- and endpoints.

We shall be concerned below with the following three well-known
models, in which no rigorous answer has previously been given to the
question of whether or not accessible paths exist with high probability.

\begin{model}[{[Unconstrained house of cards (HoC)]}]\label{mo1}
This model is originally attributed to Kingman \cite{Ki}. In the form
we consider below, it was first studied by Kauffman and Levin \cite
{KL}. We set $f(\bmm{v}^1):= 1$ and, for every other node $v \in
\mathbb{Q}_n$, independently let $f(v) \sim\operatorname{U}(0, 1)$,
the uniform distribution on the interval $[0, 1]$.
\end{model}

\begin{model}[{[Constrained house of cards (CHoC)]}]\label{mo2}
This variant seems to have been considered only more recently; see, for
example, \cite{Klo} and \cite{CH}. The only difference from
Model \ref{mo1} is that we fix $f(\bmm{v}^0):= 0$.
\end{model}

\begin{model}[{[Rough Mount Fuji (RMF)]}]\label{mo3}
This model was first proposed in \cite{A}; see also \cite{FWK}. 
For each $v \in\mathbb{Q}_n$, one lets
%
\begin{equation}
\label{eq:defrmf} f(v) = \theta\cdot d\bigl(v, \bmm{v}^0\bigr) + \eta(v),
\end{equation}
where $\theta=\theta_n$ is a positive number called the \emph{drift},
 $d(\cdot, \cdot)$ denotes Hamming distance and the $\eta(v)$ are
independent random variables of some fixed distribution. In other
words, one first assigns a fitness to each node at random, according to
$\eta$ and independent of all other nodes. Then the fitness of each
node is shifted upward by a fixed multiple of the Hamming distance from~$\bmm{v}^0$.
\end{model}

Before proceeding, it is worth noting that the above models are also of
interest in physics in the context of so-called \emph{spin glasses}
\cite{MPV}. In this setting, each node of $\mathbb{Q}_n$ represents a
point in the state space of all possible configurations of spins in a
disordered magnet. The analogue of fitness is in this case energy, or
more precisely ``energy times $-1$.'' Accessible paths (not necessarily
from $\bmm{v}^0$ to $\bmm{v}^1$) correspond to trajectories in which
energy decreases monotonically, and which are therefore easily
accessible even at zero temperature. The HoC model appears in the spin
glass context as Derrida's random energy model (REM), and the RMF-model
is a REM in an external magnetic field. For further discussion of the
connection between fitness landscapes and spin glasses, see \cite{FK}.

In all three models, the basic random variable of interest is the
number $X = X(n)$ of accessible paths. One thinks of $\bmm{v}^0$ as the
starting point of some evolutionary process, and $\bmm{v}^1$ as the
desirable endpoint. The HoC model is often referred to as a ``null
model'' for evolution, since the fitnesses of all nodes other than
$\bmm
{v}^1$ are assigned at random and independently of one another. No
mechanism is prescribed which might push an evolutionary process in any
particular direction. The CHoC model is not much better, though it does
specify that the starting point is a global fitness minimum. The RMF
model is a very natural and simple way to introduce an ``arrow of
evolution,'' since the drift factor implies that successive $0
\rightarrow1$ mutations will tend to increase fitness.

It seems intuitively obvious that the number $X$ of accessible
paths should, on average, be much higher in RMF than in HoC, with the
CHoC model lying somewhere in between. One should be a little careful
here, since in RMF, the node $\bmm{v}^1$ is not assumed to be a global
fitness maximum. Nevertheless, it is easy to verify that $\mathbb{E}[X]
= 1$ in HoC,
$\mathbb{E}[X]= n$ in CHoC, whereas in many situations $\mathbb{E}[X]$
grows super-exponentially with $n$ in RMF; see \cite{FKVK}, along with
Propositions \ref{prop:EX} and \ref{vertexremovalprop} below. Of more
interest, however, is the quantity $P = P(n)$, which is the probability
of there being at least one accessible path, that is, $P = \mathbb{P}(X
> 0)$. The idea here is that, as long as \emph{some} accessible path
exists, then evolution will eventually find it. The quantity $P$ has
been simulated in the biophysics literature. In \cite{FKVK} it was
conjectured explicitly that $P \rightarrow0$ in the HoC model, and
that $P \rightarrow1$ in the RMF model, when $\eta$ is a normal
distribution and $\theta$ is any positive constant. In \cite{CH}, the
CHoC model was simulated for $n \leq13$, and the authors conjecture,
if somewhat implicitly, that $P$ is monotonic decreasing in $n$ and
approaches a limiting value close to $0.7$. In \cite{FKVK}, simulations
were continued up to $n = 19$, and these indicated clearly that $P$ was
not, after all, monotonic decreasing. The authors abstain from making
any explicit conjecture about the limiting behavior of $P$ in CHoC.

Our main results below resolve all these issues. A crucial idea is
to consider the following slight generalization of the CHoC model:

\begin{model}[{[$\alpha$-Constrained House of Cards ($\alpha$-HoC)]}]\label{mo4}
Let $\alpha\in[0, 1]$. In this model, fitnesses are assigned as in
the CHoC model, with the exception that we set $f(\bmm{v}^0):= \alpha
$. Hence, CHoC is the case $\alpha= 0$.
\end{model}

For $\alpha\in[0, 1]$, let $P(n, \alpha)$ denote the probability of
there being an accessible path in the $\alpha$-HoC model. To simplify
notation below, we define $P(n, \alpha)=P(n, 0)$ for $\alpha<0$ and
$P(n, \alpha)=P(n, 1)$ for $\alpha>1$. Note that $P(n, \alpha)$
decreases as $\alpha$ increases. Our first main result is the following:

\begin{theorem}\label{thm:aconstr}
Let $\varepsilon= \varepsilon_n >0$. If $n\varepsilon_n \rightarrow
\infty$, then
%
\begin{equation}
\label{eq:aconstr1} \lim_{n\rightarrow\infty} P \biggl(n, \frac{\ln n}{n} -
\varepsilon _n \biggr) = 1
\end{equation}
and
%
\begin{equation}
\label{eq:aconstr2} \lim_{n\rightarrow\infty} P \biggl(n, \frac{\ln n}{n} +
\varepsilon _n \biggr) = 0.
\end{equation}
\end{theorem}



It follows immediately that $P \rightarrow1$ in the CHoC model and
that $P(n, \alpha)\rightarrow0$ for any strictly positive constant
$\alpha$. The above result says a lot more, however. It shows that
there is a very sharp threshold at $\alpha= \alpha_n = \frac{\ln
n}{n}$ for the existence of accessible paths in the $\alpha$-HoC model.
Theorem \ref{thm:aconstr} will be proven in Section~\ref{section:proofaconstr}. We have the following immediate corollary for HoC:

\begin{corollary}
Let $X$ denote the number of accessible paths
in the HoC model. Then
%
\begin{equation}
\mathbb{P} ( X>0 ) \sim\frac{\ln n}{n}.
\end{equation}
\end{corollary}
\begin{pf}
As $P(n, \alpha)$ is decreasing in $\alpha$ we know that, for any
$\alpha\in[0, 1]$, $\mathbb{P} ( X>0  ) \geq\alpha P(n,
\alpha)$. Picking $\alpha= \frac{\ln n}{n}-\varepsilon_n$ where
$n\varepsilon_n$ tends to infinity sufficiently slowly, it follows from
Theorem \ref{thm:aconstr} that $\mathbb{P} ( X>0  )
\gtrsim
\frac{\ln n}{n}$.

To get the upper bound, let $\alpha= \frac{\ln n}{n}$. Now if the
hypercube has accessible paths, then either $\bmm{v}^0$ has fitness at
most $\alpha$, or there is an accessible path where all nodes involved
have fitness at least $\alpha$. Obviously the former event occurs with
probability~$\alpha$. Concerning the latter, if
%
\begin{equation}
\label{path} \bmm{v}^0 \rightarrow v_1 \rightarrow\cdots
\rightarrow v_{n-1} \rightarrow\bmm{v}^1
\end{equation}
is any path, then the probability of all nodes along it having fitness
at least $\alpha$ is $(1-\alpha)^{n}$. The probability of fitness being
increasing along the path is $1/n!$. Since there are $n!$ possible
paths of the form \eqref{path}, it follows from a union bound that
%
\begin{equation}
\mathbb{P}(X > 0) \leq\alpha+ n!\frac{(1-\alpha)^n}{n!} \leq\frac
{\ln n}{n} +
\frac{1}{n}.
\end{equation}
\upqed\end{pf}

Another corollary of Theorem \ref{thm:aconstr} concerns the
distribution of the number of accessible paths in $\alpha$-HoC for
$\alpha=\frac{\ln n}{n} - \epsilon_n$, where $n\epsilon
_n\rightarrow
\infty$. It is straightforward to show that the expected number of
paths in $\alpha$-HoC is $n(1-\alpha)^{n-1}$ (see Proposition~\ref
{prop:EX}), which, for this choice of $\alpha$, is $\sim e^{n\epsilon
_n}$. We have the following result:
%
\begin{corollary}\label{cor:aHoCconc}
Let $X$ denote the number of accessible paths in $\alpha$-HoC for
$\alpha= \frac{\ln n}{n} - \epsilon_n$ where $n\epsilon_n
\rightarrow
\infty$. If $w_n\rightarrow\infty$, then
%
\begin{equation}
\lim_{n\rightarrow\infty}\mathbb{P} \biggl(\frac{1}{w_n}\mathbb{E}[X]
\leq X \leq w_n \mathbb{E}[X] \biggr) = 1.
\end{equation}
%
\end{corollary}
Corollary \ref{cor:aHoCconc} will be proven in Section~\ref{subsection:aHoCconc}.

Our second main result concerns the RMF model. For any function
$f\dvtx \mathbb{R} \rightarrow\mathbb{R}$, recall that the \emph{support}
of $f$, denoted $\operatorname{Supp}(f)$, is the set of points at which $f$ is
nonzero,{\footnote{Sometimes in the mathematical literature, the
support of a function is defined to be the closure of this set.}} that
is, $\operatorname{Supp}(f) = \{x \dvtx f(x) \neq0\}$. We say that $f$ has \emph{connected
support} if $\operatorname{Supp}(f)$ is a connected subset of $\mathbb{R}$. Our result
is the following:

\begin{theorem}\label{thm:rmf}
Let $\eta$ be any probability distribution whose p.d.f. is continuous
on its support and whose support is connected. Let $\theta_n$ be any
strictly positive function of $n$ such that $n \theta_n \rightarrow
\infty$ as $n \rightarrow\infty$. Then in the model \eqref{eq:defrmf},
$P(n)$ tends to one as $n \rightarrow\infty$.
\end{theorem}

This result is proven in Section~\ref{section:rmf}. The proof follows
similar lines to that of Theorem \ref{thm:aconstr}, but the analysis is
somewhat simpler.

\begin{remark}\label{re1.6}
More generally, the proof of Theorem \ref{thm:rmf} presented in this
article holds for any distribution $\eta$ that satisfies, with notation
taken from Section~\ref{section:rmf}, $\kappa_{\eta,\delta} = \inf_{I
\subseteq I_\delta} \frac{1}{l(I)} \int_I \eta(x) \,dx > 0$ for any
$\delta\in(0, 1)$. This condition essentially states that $\eta$ is
not allowed to have ``isolated modes.'' For instance, it is satisfied
for any unimodal distribution.
\end{remark}

\section{Results for the HoC models}\label{section:proofaconstr}
For each path $i$ from $\bmm{v}^0$ to $\bmm{v}^1$ let $X_i$ be the
indicator function of the event that $i$ is accessible, and let $X=\sum_i X_i$ denote the number of accessible paths from $\bmm{v}^0$ to
$\bmm
{v}^1$. Furthermore, given a path $i$ from $\bmm{v}^0$ to $\bmm{v}^1$
in the $n$-dimensional hypercube, let $T(n, k)$ denote the number of
paths from $\bmm{v}^0$ to $\bmm{v}^1$ that intersect $i$ in exactly
$k-1$ interior nodes (by symmetry, this is independent of $i$).

\begin{proposition}\label{prop:EX}
Let $X$ denote the number of accessible paths in the $\alpha$-HoC
model. Then
%
\begin{equation}
\mathbb{E}[X] = n(1-\alpha)^{n-1}.
\end{equation}
\end{proposition}
\begin{pf}
There are $n!$ paths through the hypercube. A path is accessible if all
$n-1$ interior nodes have fitness at least $\alpha$, and the fitness of
the interior nodes is increasing along the path. This occurs with
probability $(1-\alpha)^{n-1}/(n-1)!$.
\end{pf}

Note that for $\alpha= \frac{\ln n}{n} + \epsilon_n$, the proposition
implies that the expected number of accessible paths tends to $0$ for
any sequence $\epsilon_n$ satisfying $n\epsilon_n \rightarrow\infty$.
This directly implies equation \eqref{eq:aconstr2}. Similarly, for
$\alpha= \frac{\ln n}{n}-\epsilon_n$ where $n\epsilon_n\rightarrow
\infty$, the expected number of paths tends to infinity.

To show the remaining part of Theorem \ref{thm:aconstr}, that the
probability of there being at least one accessible path tends to $1$ in
the case $\alpha= \frac{\ln n}{n} - \epsilon_n$, we will begin by
showing that the probability is at least $\frac{1}{4}-o(1)$ by the
second moment method. In Section~\ref{subsection:bootstrap} we will
then provide a proof that the probability must tend to 1.

%
\begin{lemma}\label{lemma:alon}
Let $X$ be a random variable with finite expected value and finite and
nonzero second moment. Then
%
\begin{equation}
\mathbb{P} ( X \neq0 ) \geq\frac{\mathbb
{E}[X]^2}{\mathbb{E}[X^2]}.
\end{equation}
\end{lemma}
\begin{pf}
Let $1_{X\neq0}$ denote the indicator function of $X\neq0$. Then, by
the Cauchy--Schwarz inequality, $\mathbb{E}[X]^2 = \mathbb
{E}[1_{X\neq
0}X]^2 \leq\mathbb{E}[1_{X\neq0}^2] \cdot\mathbb{E}[X^2] = \mathbb
{P}(X\neq0) \cdot\mathbb{E}[X^2]$.
\end{pf}
See also Exercise 4.8.1 in \cite{AS}.

\begin{proposition}\label{prop:EX2ub}
Let $i$ and $j$ be paths with exactly $k-1$ interior nodes in common. Then
%
\begin{equation}
\label{eq:XiXj} \mathbb{E} [ X_i X_j ] \leq
\frac{ {2n-2k \choose n-k}
(1-\alpha)^{2n-k-1} }{(2n-k-1)!},
\end{equation}
where equality holds if the nodes where $i$ and $j$ differ are
consecutive along the paths, that is, if $i$ and $j$ diverge at most
once. Furthermore,
%
\begin{equation}
\label{eq:theonewiththesum} \mathbb{E}\bigl[X^2\bigr] \leq\sum
_{k=1}^n n! T(n, k) \frac{ {2n-2k \choose n-k}
(1-\alpha)^{2n-k-1} }{(2n-k-1)!}.
\end{equation}
\end{proposition}
\begin{pf}
The event that $i$ and $j$ are both accessible occurs if all $2n-k-1$
interior nodes have fitness at least $\alpha$ and the fitnesses of the
interior nodes are ordered in such a way that fitness increases along
both paths.

Conditioned on the event that all interior nodes have fitness at least
$\alpha$, all possible ways in which the fitnesses of the interior
nodes can be ordered are equally likely. This implies that the
probability that both paths are accessible is $(1-\alpha
)^{2n-k-1}/(2n-k-1)!$ times the number of ways to order the fitnesses
of the interior nodes such that fitness increases along both paths.

To count the number of ways this can be done we color the numbers $1,
\ldots, 2n-k-1$ in the following way: The number $l$ is colored gray if
the interior node with the $l$th smallest fitness is contained in both
paths, red if it is only contained in $i$ and blue if only in $j$. Note
that $i$ and $j$ uniquely determine which numbers must be gray for a
valid order, and that any coloring corresponds to at most one order.

Clearly, any coloring corresponding to a valid order colors half of the
nongray numbers red and half blue, which implies that there can be at
most ${2n-2k \choose n-k}$ such orders. Furthermore, if $i$ and $j$
diverge at most once, one can always construct a valid order from such
a coloring, so in this case there are exactly ${2n-2k \choose n-k}$
such orders.

As the number of ordered pairs of paths that intersect in exactly $k-1$
interior nodes is $n!T(n, k)$, \eqref{eq:theonewiththesum} follows from
this estimate.
\end{pf}

\subsection{Useful formulas for $T(n,k)$}


The numbers $T(n, k)$ already appear in the mathematical literature.
The usual terminology is that $T(n, k)$ is the number of permutations
of $\{1, 2, \ldots, n\}$ with $k$ \emph{components}, where the number of
components of a permutation $\pi_1 \pi_2 \cdots\pi_n$ is defined as
the number of choices for $1\leq s \leq n$ such that $\pi_1 \pi_2
\cdots\pi_s$ is a permutation of $\{1, 2, \ldots, s\}$. In terms of
paths in $\Q_n$, we can represent each path from $\bmm{v}^0$ to $\bmm
{v}^1$ by a permutation $\pi_1 \pi_2 \cdots\pi_n$ of $\{1, 2, \ldots,
n\}$ where $\pi_s$ denotes which coordinate to increase in step $s$. If
we let $i$ be the path represented by the identity permutation, then a
path $j$, represented by $\pi_1 \pi_2 \cdots\pi_n$, intersects $i$ in
step $s \geq1$\vadjust{\goodbreak} if and only if $\pi_1 \pi_2 \cdots\pi_s$ is a
permutation of $\{1, 2, \ldots, s\}$. This means that, if $\pi_1 \pi_2
\cdots\pi_n$ has $k$ components, then $i$ and $j$ intersect in $k-1$
interior nodes (the $k$th component corresponds to $s=n$). We can thus
consider a component as an interval $[s, t]$ where $i$ and $j$
intersect in steps $s$ and $t$, but at no step in between.

An alternative formulation is that $T(n, k)$ is the number of
permutations of $\{1, 2,\ldots, n\}$ with $k-1$ \emph{global descents}. A
global descent in a permutation $\pi_1 \pi_2 \cdots\pi_n$ of $\{1,
2,\ldots, n\}$ is a number $t \in[1, n-1]$ such that $\pi_i > \pi_j$
for all $i \leq t$ and $j > t$. There is a simple 1--1 correspondence
between permutations with $k$ components and those with $k-1$ global
descents obtained by reading a permutation backward. In other words,
$\pi_1 \pi_2 \cdots\pi_n$ has $k-1$ global descents if and only if
$\pi
_n \pi_{n-1} \cdots\pi_1$ has $k$ components.

There is a database of the numbers $T(n, k)$ for small $n$ and $k$; see
\cite{O2}.
Comtet's book \cite{Co2} contains a couple of exercises and an implicit
recursion formula for $T(n, k)$. Comtet has also performed a detailed
asymptotic analysis of the numbers $T(n, 1)$ in \cite{Co1}.
Permutations with one component (i.e., no global descents) are
variously referred to as \emph{connected, indecomposable, irreducible}.
These seem to crop up quite a lot; see \cite{O1}.
However, estimates of the numbers $T(n, k)$ for general $n$ and $k$
like those in Propositions \ref{prop:boundTnksmallk} and \ref
{prop:boundTnklargek} below do not appear to have been obtained before.

\begin{proposition}\label{prop:Tn1}
$T(n, 1)$ is uniquely defined by
%
\begin{equation}
n! = \sum_{k=1}^n T(k, 1) (n-k)!.
\end{equation}
\end{proposition}
\begin{pf}
Given a path $i$ through $\mathbb{Q}_n$, the number of paths $j$ that
intersect $i$ for the first time in step $k$ is $T(k, 1)(n-k)!$. As any
path through $\mathbb{Q}_n$ intersects $i$ for the first time after
between $1$ and $n$ steps, the proposition follows.
\end{pf}

\begin{proposition}\label{prop:Tn1bound}
%
\begin{equation}
n! \biggl( 1-O \biggl(\frac{1}{n} \biggr) \biggr) \leq T(n, 1) \leq n!.
\end{equation}
\end{proposition}
\begin{pf}
By definition, $T(n, 1) \leq n!$. Using this, Proposition \ref
{prop:Tn1} implies that $T(n, 1)$ is at least $n!-\sum_{k=1}^{n-1}
k!(n-k)! = n! - O ( (n-1)!  )$.
\end{pf}

\begin{proposition}\label{prop:Tnksumovers}
%
\begin{equation}
\label{eq:Tnksumovers} T(n, k) = \mathop{\sum_{s_1,\ldots, s_k \geq1}}_{ s_1+\cdots+ s_k = n}
T(s_1, 1)\cdots T(s_k, 1).
\end{equation}
\end{proposition}
\begin{pf}
Given a path $i$, the number of paths that intersect $i$ for the first
time after $s_1$ steps, for the second time after $s_2$ more steps and
so on up to the last time (at $\bmm{v}^1$) after $n$ steps is $T(s_1,
1)\cdots T(s_{k-1}, 1) \cdot T(n-s_1-\cdots-s_{k-1}, 1)$. Let
$s_k=n-s_1-\cdots-s_{k-1}$. $T(n, k)$ is obtained by summing over all
possible values of $s_1,\ldots, s_k$.
\end{pf}

\begin{proposition}\label{prop:Tnk}
For $k\geq2$, $T(n, k)$ satisfies
%
\begin{equation}
T(n, k) = \sum_{s=1}^{n-k+1} T(s, 1) T(n-s,
k-1).
\end{equation}
\end{proposition}
\begin{pf}
It follows by induction that this sum equals the right-hand side in~\eqref{eq:Tnksumovers}.
\end{pf}

\subsection{Upper bounds for $T(n,k)$}

\begin{proposition}\label{prop:grouponlargest}
For any $n\geq k \geq1$,
%
\begin{equation}
\label{eq:grouponlargest} T(n, k) \leq k \sum \Biggl(\Biggl(n-\sum
_{j=1}^{k-1} s_j\Biggr)! \prod
_{j=1}^{k-1} s_j! \Biggr),
\end{equation}
where the first sum goes over all $(k-1)$-tuples of integers $s_1,
\ldots, s_{k-1}$ such that $s_j \geq1$ for all $j$ and $\max_j s_j \leq
n-\sum_j s_j$.
\end{proposition}
\begin{pf}
Consider the formula for $T(n, k)$ in Proposition \ref
{prop:Tnksumovers}. By symmetry, $T(n, k)$ is at most $k$ times the
contribution from terms where $s_j\leq s_k$ for $j=1, \ldots, k-1$. The
proposition follows by applying $T(s, 1) \leq s!$.
\end{pf}

\begin{proposition}\label{prop:boundTnksmallk}
There is a positive constant $c$ such that for all $n \geq k \geq1$,
%
\begin{equation}
T(n, k) \leq k (n-k+1)! e^{c(k-1)/(n-k+1)}.
\end{equation}
\end{proposition}
\begin{pf}
We use Proposition \ref{prop:grouponlargest} and make the following
approximations:
\begin{itemize}
\item Substitute $(n-\sum_j s_j)!$ by $\beta^{n-\sum_j s_j}$ where
$\beta=  ( (n-k+1)!  )^{1/(n-k+1)}$. It follows from
log-convexity of $l!$ that $\beta^{l} \geq l!$ for any $0 \leq l \leq n-k+1$.
\item Let all $s_j$ go from $1$ to $\lfloor(n-k+1)/2+1 \rfloor$.
\end{itemize}
This yields
%
\begin{equation}
\label{eq:Tnkbound} T(n, k) \leq k (n-k+1)! \Biggl(\sum_{s=1}^{\lfloor(n-k+1)/2+1 \rfloor}
s! \beta^{1-s} \Biggr)^{k-1}.
\end{equation}

We now claim that the sum in the above expression is always less than
$1+c/(n-k+1)$ for sufficiently large $c$. Indeed,
\begin{eqnarray*}
&&\sum_{s=1}^{\lfloor(n-k+1)/2+1 \rfloor} s! \beta^{1-s} \\
&&\qquad=
1 + 2\beta ^{-1} + \beta^{-1}\sum
_{t=1}^{\lfloor(n-k+1)/2-1 \rfloor} t! (t+1) (t+2) \beta^{-t}
\\
&&\qquad\leq 1 + 2\beta^{-1}
\\
&&\qquad\quad{} + e\beta^{-1} \sum_{t=1}^{\lfloor(n-k+1)/2-1 \rfloor}
\sqrt{t} (t+1) (t+2) \biggl( \frac{n-k+1}{2e} \biggr)^t \biggl(
\frac
{n-k+1}{e} \biggr)^{-t}
\\
&&\qquad\leq 1 + 2\beta^{-1} + e\beta^{-1} \sum
_{t=1}^{\infty} \sqrt{t} (t+1) (t+2)2^{-t}
\\
&&\qquad\leq 1 + c(n-k+1)^{-1}.
\end{eqnarray*}
Here we have used that $(n-k+1)/e \leq\beta\leq(n-k+1)$ and that $n!
\leq e n^{n+1/2}e^{-n}$, which follows from standard estimates of
factorials. 

The proposition now follows from this result together with \eqref{eq:Tnkbound}.
\end{pf}

\begin{proposition}\label{prop:Tnkfixedl}
For any fixed $l$ there is a constant $C_l>0$ such that
%
\begin{equation}
T(n, n-l) \leq C_l n^l
\end{equation}
for all $n\geq1$.
\end{proposition}
\begin{pf}
We may, without loss of generality, assume that $n \geq2l$.

Recall the formula for $T(n, n-l)$ in Proposition \ref
{prop:Tnksumovers}. As $s_1, \ldots, s_{n-l} \geq1$ and $s_1+\cdots
+s_{n-l}=n$ it is easy to see that all but at most $l$ variables are
equal to $1$. This implies that $T(n, n-l)$ is at most ${n-l \choose
l}$ times the contribution from all terms where $s_{l+1}=\cdots=
s_{n-l}=1$. Using $T(1, 1)=1$, we get
%
\begin{equation}
T(n, n-l) \leq\pmatrix{n-l \cr l} \mathop{\sum_{s_1,\ldots, s_l \geq1}}_
{ s_1+\cdots+ s_l = 2l}
T(s_1, 1)\cdots T(s_l, 1) \leq
C_l n^l.
\end{equation}
\upqed\end{pf}

\begin{proposition}\label{prop:boundTnklargek}
For sufficiently large $c$, we have
%
\begin{equation}
T(n, n-l) \leq c(l+1) \biggl(\frac{n+2l}{5} \biggr)^l.
\end{equation}
\end{proposition}

\begin{pf}
Let
%
\begin{equation}
S(n, n-l) = (l+1) \biggl( \frac{n+2l}{5} \biggr)^l,
\end{equation}
that is,
%
\begin{equation}
S(n, k) = (n-k+1) \biggl( \frac{3n-2k}{5} \biggr)^{n-k}.
\end{equation}

We will begin by showing that $S(n, k)$ satisfies
%
\begin{equation}
\label{eq:Tnklargekrec} S(n, k) \geq\sum_{i=1}^{n-k+1}
i! S(n-i, k-1)
\end{equation}
for $k>1$ and sufficiently large $n-k$. Here we have
\begin{eqnarray*}
&&\sum_{i=1}^{n-k+1} i! S(n-i, k-1)\\
&&\qquad= \sum
_{i=1}^{n-k+1} i! (n-k+2-i) \biggl(
\frac{3n-2k-3i+2}{5} \biggr)^{n-k-i+1}
\\
&&\qquad\leq(n-k+1) \biggl(\frac{3n-2k-1}{5} \biggr)^{n-k}
\\
&&\qquad\quad{} + \sum_{i=2}^{n-k+1} i! (n-k+1) \biggl(
\frac{3n-2k}{5} \biggr)^{n-k-i+1}
\\
&&\qquad= S(n, k) \Biggl( \biggl( 1 - \frac{1}{3n-2k} \biggr)^{n-k} + \sum
_{i=2}^{n-k+1} i! \biggl(\frac{3n-2k}{5}
\biggr)^{-i+1} \Biggr),
\end{eqnarray*}
where
\begin{eqnarray*}
\biggl( 1 - \frac{1}{3n-2k} \biggr)^{n-k} &\leq&\exp \biggl( -
\frac
{n-k}{3n-2k} \biggr)
\\
&\leq&\exp \biggl( - \frac{n-k}{3n} \biggr)
\leq\max \biggl(\frac{1}{2}, 1-\frac{n-k}{6n} \biggr)
\end{eqnarray*}
and
\begin{eqnarray*}
&&\sum_{i=2}^{n-k+1} i! \biggl(
\frac{3n-2k}{5} \biggr)^{-i+1}
\\
&&\qquad \leq\frac{10}{3n-2k} + \frac{5}{3n-2k}\sum_{j=1}^{n-k-1}
j!(j+1) (j+2) \biggl(\frac{3n-2k}{5} \biggr)^{-j}
\\
& &\qquad\leq\frac{10}{3n-2k} + \frac{5e}{3n-2k}\sum_{j=1}^{\infty}
\sqrt {j}(j+1) (j+2) \biggl( \frac{n-k}{e} \biggr)^j \biggl(
\frac
{3n-2k}{5} \biggr)^{-j}
\\
& &\qquad\leq\frac{1}{n} \Biggl( 10+5e \sum_{j=1}^\infty
\sqrt{j}(j+1) (j+2) \biggl(\frac{5}{3e} \biggr)^j \Biggr)
\\
& &\qquad= \frac{C}{n}.
\end{eqnarray*}
It follows directly that~\eqref{eq:Tnklargekrec} holds for $k>1$ and
$n-k \geq6C$.

Now, if we can choose $c$ so that $T(n, k) \leq cS(n, k)$ for $k=1$ and
for $n-k<6C$, the proposition will follow from Proposition \ref
{prop:Tnk} by induction on $k$. Hence it suffices to show the
proposition for these two cases.

For $k=1$, the inequality holds for sufficiently large $c$ by the fact that
\begin{eqnarray*}
\frac{T(n, 1)}{S(n, 1)} &\leq&\frac{n!}{n  ({(3n-2)}/{5}
)^{n-1} }
\\
&\leq& e \sqrt{n} \biggl(\frac{n}{e} \biggr)^n
\frac{1}{n
(
{(3n-2)}/{5} )^{n-1} }
\\
&= &\frac{3e}{5} \sqrt{n} \biggl( \frac{5}{3e} \biggr)^n
\biggl( 1 - \frac{2}{3n} \biggr)^{-n+1}
\\
&\rightarrow&0\qquad \mbox{as } n\rightarrow\infty.
\end{eqnarray*}
For $n-k<6C$, just apply Proposition \ref{prop:Tnkfixedl}.
\end{pf}

\subsection{Computing $\mathbb{E}[X^2]$}
Pick $\delta> 0$ sufficiently small. We divide the sum in~\eqref
{eq:theonewiththesum} into the contribution from $k \leq(1-\delta)n$
and that from $k > (1-\delta)n$:
%
\begin{eqnarray}
\label{eq:dividesum} &&\sum_{k=1}^{n} n! T(n, k)
\frac{ {2n-2k \choose n-k} (1-\alpha
)^{2n-k-1} }{ (2n-k-1)! }
\nonumber\\
&&\qquad = \sum_{k=1}^{(1-\delta)n} n! T(n, k)
\frac{ {2n-2k \choose n-k}
(1-\alpha)^{2n-k-1} }{ (2n-k-1)! }
\nonumber
\\[-8pt]
\\[-8pt]
\nonumber
&&\qquad\quad{} + \sum_{l=0}^{\delta n} n! T(n,
n-l) \frac{ {2l \choose l} (1-\alpha)^{n+l-1} }{ (n+l-1)! }
\\
&&\qquad:= S_1 + S_2.\nonumber
\end{eqnarray}

\begin{proposition}\label{prop:Tnkconstantk}
For $k$ constant and $\alpha=o(1)$
%
\begin{equation}
n! T(n, k) \frac{ {2n-2k \choose n-k} (1-\alpha)^{2n-k-1} }{ (2n-k-1)!
} \sim k 2^{1-k} n^2 (1-
\alpha)^{2n}.
\end{equation}
\end{proposition}
\begin{pf}
A simple lower bound on $T(n, k)$ is the number of permutations with
$k$ components where all but one component contains exactly one
element. For sufficiently large $n$ this is given by $k T(n-k+1, 1)$,
which by Proposition \ref{prop:Tn1bound} is $\sim k(n-k+1)!$.
Furthermore, from Proposition \ref{prop:boundTnksmallk} we know that
$T(n, k)$ is most $ (1+o(1) )k(n-k+1)!$. Hence for constant
$k$, $T(n, k) \sim k(n-k+1)!$. The proposition now follows from
standard estimates of factorials.
\end{pf}

\begin{proposition}\label{prop:divsumsmallk}
Let $\alpha= o(1)$. For any $0 < \delta< 1$, we have $S_1 \sim4 n^2
(1-\alpha)^{2n}$.
\end{proposition}
\begin{pf}
From Proposition \ref{prop:boundTnksmallk} it follows that there is a
constant $C_\delta$ such that
$T(n, k) \leq C_\delta k(n-k+1)!$ whenever $k \leq(1-\delta)n$. Using
this we have
%
\begin{equation}
n! T(n, k) \frac{ {2n-2k \choose n-k} }{ (2n-k-1)! } \leq C_\delta n! k (n-k+1)!
\frac{ {2n-2k \choose n-k} }{ (2n-k-1)! }
\end{equation}
for all $k \leq(1-\delta)n$. Now by extensive use of Stirling's
formula there is a constant $C>0$ such that
\begin{eqnarray*}
&&C_\delta n! k (n-k+1)! \frac{ {2n-2k \choose n-k} }{ (2n-k-1)! }
\\
& &\qquad\leq C_\delta C k \sqrt{n} \biggl( \frac{n}{e}
\biggr)^n \sqrt{n-k} \biggl( \frac{n-k}{e} \biggr)^{n-k}(n-k+1)\\
&&\qquad\quad{}\times
\frac{ ({4^{n-k} }/{
\sqrt {n-k} })(2n-k) }{ \sqrt{2n-k} ( {(2n-k)}/{e}  )^{2n-k} }
\\
&&\qquad = C_\delta C k (n-k+1)\sqrt{n(2n-k)} 2^{-k}\\
&&\qquad\quad{}\times \biggl(
\biggl(1-\frac
{k}{n} \biggr)^{{n}/{k}-1} \biggl( 1-\frac{k}{2n}
\biggr)^{-{2n}/{k}+1} \biggr)^k,
\end{eqnarray*}
where
\begin{eqnarray*}
\biggl(1-\frac{k}{n} \biggr)^{{n}/{k}-1} \biggl( 1-\frac
{k}{2n}
\biggr)^{-{2n}/{k}+1} &\leq &\biggl(1-\frac{k}{2n} \biggr)^{{2n}/{k}-2}
\biggl( 1-\frac{k}{2n} \biggr)^{-{2n}/{k}+1}
\\
&= &\biggl( 1-\frac{k}{2n} \biggr)^{-1}
\\
&\leq& \biggl( 1-\frac{1-\delta}{2} \biggr)^{-1}
\\
&=& \frac{2}{1+\delta}.
\end{eqnarray*}
This means that, for all $\delta>0$, there exists a constant $C_\delta
'$ such
that, for $k\leq(1-\delta)n$ and sufficiently large $n$, we have
%
\begin{eqnarray}
&&n! T(n, k) \frac{ {2n-2k \choose n-k} (1-\alpha)^{2n-k-1} }{ (2n-k-1)!
}
\nonumber
\\[-8pt]
\\[-8pt]
\nonumber
&&\qquad\leq C_\delta'
n^2(1-\alpha)^{2n} k (1+\delta )^{-k}(1-
\alpha)^{-k}.
\end{eqnarray}

Since $\sum k(1+\delta)^{-k}(1-\alpha)^{-k}$ converges for sufficiently
small $\alpha$ we have shown that $S_1=O (n^2(1-\alpha
)^{2n}
)$. Furthermore, if we assume that $n$ is sufficiently large so that
$(1+\delta)(1-\alpha)\geq(1+\frac{\delta}{2})$, then as the terms in
the sum
%
\begin{equation}
\sum_{k=1}^{(1-\delta)n}\frac{1}{n^2(1-\alpha)^{2n} } n! T(n,
k) \frac{
{2n-2k \choose n-k} (1-\alpha)^{2n-k-1} }{ (2n-k-1)! }
\end{equation}
are dominated by the terms in
%
\begin{equation}
\sum_{k=1}^\infty C_\delta'
k \biggl(1+\frac{\delta}{2} \biggr)^{-k},
\end{equation}
which converges, it follows by dominated convergence together with
Proposition~\ref{prop:Tnkconstantk} that
\begin{eqnarray}
\sum_{k=1}^{(1-\delta)n}\frac{1}{n^2(1-\alpha)^{2n} } n! T(n,
k) \frac{
{2n-2k \choose n-k}(1-\alpha)^{2n-k-1} }{ (2n-k-1)! }\longrightarrow \sum_{k=1}^\infty
k 2^{1-k} = 4\nonumber\\
 \eqntext{\mbox{as } n\rightarrow\infty.\qquad}
\end{eqnarray}
\upqed\end{pf}

\begin{proposition}\label{prop:divsumlargek}
For sufficiently small $\delta> 0$ and $\alpha=o(1)$, we have $S_2 =
O (n (1-\alpha)^n  )$.
\end{proposition}
\begin{pf}
Using Proposition \ref{prop:boundTnklargek} there is a constant $C$
such that this sum is bounded by
\begin{eqnarray*}
\sum_{l=0}^{\delta n} n! T(n, n-l)
\frac{ {2l \choose l} (1-\alpha
)^{n+l-1} }{ (n+l-1)! } &\leq& C\sum_{l=0}^{\delta n} n!
(l+1) \biggl( \frac{n+2l}{5} \biggr)^{l} \frac{ {2l \choose l} (1-\alpha)^{n+l-1} }{
(n+l-1)! }
\\
&\leq &C (1-\alpha)^{n-1} \sum_{l=0}^{\delta n}
n^{1-l}(l+1) \biggl( \frac{n+2l}{5} \biggr)^l
4^l
\\
&\leq & C n (1-\alpha)^{n-1} \sum_{l=0}^{\infty}
(l+1) \biggl( \frac
{4(1+2\delta)}{5} \biggr)^l,
\end{eqnarray*}
where the last sum clearly converges for sufficiently small $\delta$.
\end{pf}

\begin{proposition}\label{prop:EX2}
Let $X$ be the number of accessible paths in the $\alpha$-HoC model
where $\alpha= \frac{\ln n}{n}-\epsilon_n$ where $n\epsilon
_n\rightarrow\infty$. Then
%
\begin{equation}
\mathbb{E}\bigl[X^2\bigr] \sim4n^2(1-
\alpha)^{2n}.
\end{equation}
\end{proposition}
\begin{pf}
From Proposition \ref{prop:EX2ub} together with Propositions \ref
{prop:divsumsmallk} and \ref{prop:divsumlargek} we know that
%
\begin{equation}
\mathbb{E}\bigl[X^2\bigr] \leq \bigl(4 + o(1) \bigr)
n^2(1-\alpha)^{2n} + O \bigl( n (1-\alpha)^{n}
\bigr),
\end{equation}
where one can show that $n (1-\alpha)^{n} = o ( n^2 (1-\alpha)^{2n}
)$, provided $n\epsilon_n \rightarrow\infty$.

To derive a tight lower bound for $\mathbb{E}[X^2]$, consider the sum
of $\mathbb{E}[X_iX_j]$ over all pairs of paths whose number of common
interior nodes, $k-1$, is at most $\frac{n}{2}-1$ and that diverge at
most once. Expressed in terms of components of permutations, for a
fixed $i$ and $k$, the number of paths $j$ that satisfy this equals the
number of permutations with $k$ components, where all but one component
contains exactly one element. This can clearly be done in $kT(n-k+1, 1)
\sim k(n-k+1)!$ ways.

By Proposition \ref{prop:EX2ub} this yields
%
\begin{equation}
\mathbb{E}\bigl[X^2\bigr] \geq\sum_{k=1}^{n/2}
n! kT(n-k+1, 1) \frac{ {2n-2k
\choose n-k}(1-\alpha)^{2n-k-1} }{(2n-k-1)!}.
\end{equation}
Proceeding in a manner similar to the proof of Proposition \ref
{prop:divsumsmallk}, we get that
%
\begin{equation}
\sum_{k=1}^{n/2} n! kT(n-k+1, 1)
\frac{ {2n-2k \choose n-k}(1-\alpha
)^{2n-k-1} }{(2n-k-1)!} \sim4n^2(1-\alpha)^{2n},
\end{equation}
which completes the proof.\vadjust{\goodbreak}
\end{pf}

From this proof we can observe that almost all of the contributions to
$\mathbb{E}[X^2]$ come from pairs of paths we considered in the lower
bound. This implies the following:
%
\begin{corollary}\label{cor:smallcontrib}
Assume $\alpha= \frac{\ln n}{n}-\epsilon_n$ where $n\epsilon
_n\rightarrow\infty$. For any $0<\delta<1$, the contribution to
$\mathbb{E}[X^2]$ from all pairs of paths that either share more than
$(1-\delta)n$ common nodes or that diverge more than once is $o (
n^2(1-\alpha)^{2n}  )$.
\end{corollary}

\subsection{Proof of Theorem \texorpdfstring{\protect\ref{thm:aconstr}}{1.2}}\label
{subsection:bootstrap} Let $X$ as above denote the number of accessible
paths in $\alpha$-HoC, where $\alpha= \frac{\ln n}{n} - \epsilon_n$,
$0\leq\epsilon_n \leq\frac{\ln n}{n}$ and $n\epsilon_n \rightarrow
\infty$. Applying Lemma \ref{lemma:alon} to $X$ and using the
expressions for $\mathbb{E}[X]$ and $\mathbb{E}[X^2]$ from Propositions~\ref{prop:EX} and \ref{prop:EX2}, respectively, yields the lower bound
%
\begin{equation}
\label{eq:oneforth} \liminf_{n\rightarrow\infty} P \biggl(n, \frac{\ln n}{n}-
\epsilon_n \biggr) \geq\frac{1}{4}.
\end{equation}
In this subsection, we will prove that this probability can be
``bootstrapped'' up to~1, proving the remaining part of Theorem \ref
{thm:aconstr}.

\begin{lemma}\label{lemma:alphaequivalent1}
Let $0\leq a \leq1-b \leq1$, and let $f\dvtx \Q_n\rightarrow\R$ be a
fitness function whose values are generated independently according to
%
\begin{equation}
f(v) = \cases{a, &\quad $\mbox{if } v=\bmm{v}^0,$ \vspace*{2pt}
\cr
1-b, &\quad $
\mbox{if } v=\bmm{v}^1,$ \vspace*{2pt}
\cr
\sim\operatorname{U}(0, 1), &\quad$
\mbox{otherwise}$.}
\end{equation}
Then the probability of accessible paths with respect to $f$ equals
$P(n, a+b)$.
\end{lemma}
\begin{pf}
Define the function $g\dvtx \Q_n\rightarrow\R$ by setting $g(v) = f(v) + b$
if $f(v)\leq1-b$ and $g(v)=f(v)-1+b$ otherwise. Then $g(\bmm
{v}^0)=a+b$, $g(\bmm{v}^1)=1$ and $g(v) \sim\operatorname{U}(0, 1)$
independently for all other $v$, so $g$ is distributed as in $\alpha
$-HoC with $\alpha=a+b$. As this transformation only constitutes a
translation for any node on an accessible path, we see that a path is
accessible with respect to $f$ if and only if it is so with respect to $g$.
\end{pf}


\begin{proposition}\label{prop:bootstrap}
Assume there is a positive constant $C$ such that $\liminf_{n\rightarrow\infty} P(n, \frac{\ln n}{n} - \epsilon_n) \geq C$
whenever $0 \leq\epsilon_n \leq\frac{\ln n}{n}$ is a sequence
satisfying $n\epsilon_n \rightarrow\infty$. Then, the same inequality
holds if $C$ is replaced by $1-(1-\break C)(1-\frac{C}{2})$.
\end{proposition}
\begin{pf}
Let $\alpha= \frac{\ln n}{n} - \epsilon_n$. We wish to pick four
nodes, $a_1, a_2, b_1, b_2$, satisfying the following conditions:
\begin{longlist}[(iii)]
\item[(i)] $d(a_1, \bmm{v}^{0})=d(a_2, \bmm{v}^{0})=1$ and $a_1, a_2$
each has fitness in the range $[\alpha, \alpha+\epsilon_n/3]$;

\item[(ii)] $d(b_1, \bmm{v}^{1})=d(b_2, \bmm{v}^{1})=1$ and $b_1, b_2$
each has fitness at least $1-\epsilon_n/3$;

\item[(iii)] none of the four pairs $(a_i, b_j)$ are antipodal (in the
undirected hypercube).
\end{longlist}

By (i), the number of possibilities for each $a_i$ is binomially
distributed with parameters $\operatorname{Bin}(n, \epsilon_n/3)$.
Then, by (ii) and (iii), the number of options for each $b_j$ is
distributed as $\operatorname{Bin}(n-2, \epsilon_n/3)$. Since
$n\epsilon
_n/3 \rightarrow\infty$, it follows that it is possible to choose four
nodes satisfying (i)--(iii) with probability $1-o_n(1)$.

Condition on the fitness of all vertices $v$ with $d(v, \bmm{v}^0)=1$
or $d(v, \bmm{v}^1)=1$.
Let $H_1$ and $H_2$ be the induced subgraphs consisting of all nodes on
paths from $a_1$ to $b_1$ and from $a_2$ to $b_2$, respectively, and let
$H_2'$ be the induced subgraph consisting of all nodes on paths between
$a_2$ and $b_2$ that does not intersect $H_1$ in any vertex. Then $H_1$
and $H_2$ are isomorphic to $\Q_{n-2}$. Note that any accessible path
from $a_1$ to $b_1$ or $a_2$ to $b_2$ can be extended to an accessible
path from $\bmm{v}^0$ to $\bmm{v}^1$.

Let us denote the probability of accessible paths through the
respective induced subgraphs by $p_{H_1}$, $p_{H_2}$ and $p_{H_2'}$. By
construction, $H_1$ and $H_2'$ are vertex disjoint, so the events of
accessible paths through the two subgraphs are independent. By Lemma
\ref{lemma:alphaequivalent1}, $p_{H_1} = P(n-2, f(a_1)+1-f(b_1)) \geq
P(n-2, \alpha+\frac{2\epsilon_n}{3})$. It is straightforward to show
that this is still below the threshold, which implies that $p_{H_1}
\geq C - o_n(1)$.

To estimate $p_{H_2'}$, we note that a path in $H_2$ from $a_2$ to
$b_2$ is contained in $H_2'$ if and only if it ``flips the bit that is
$1$ in $a_1$ after that which is $0$ in $b_1$.'' In the cases where
there is an accessible path through $H_2$, let $\gamma$ be chosen
uniformly among all such paths. Then, by symmetry, we know that it
flips the two bits corresponding to $a_1$ and $b_1$ in the allowed
order, and is thus contained in $H_2'$, with probability $\frac{1}{2}$.
Hence $p_{H_2'} \geq\frac{1}{2} p_{H_2} = \frac{1}{2} p_{H_1}$.

As the events of accessible paths through $H_1$ and $H_2'$ are
independent, we get $P(n, \alpha) \geq1-(1-p_{H_1})(1-p_{H_2'}) -
o_n(1) \geq1-(1-C)(1-\frac{C}{2}) - o_n(1)$ and the proposition follows.
\end{pf}

Now we complete the proof of Theorem \ref{thm:aconstr}. By equation
\eqref{eq:oneforth} and repeated use of Proposition \ref
{prop:bootstrap} we can construct a sequence $\{C_k\}_{k=0}^\infty$
such that $C_k \rightarrow1$ and $\liminf_{n\rightarrow\infty}P(n,
\alpha) \geq C_k$ for all $k$. Hence we must have $\liminf_{n\rightarrow\infty} P(n, \alpha) = 1$.


\subsection{Proof of Corollary \texorpdfstring{\protect\ref{cor:aHoCconc}}{1.4}}\label
{subsection:aHoCconc}

Similarly to the proof of Theorem \ref{thm:aconstr}, that of Corollary
\ref{cor:aHoCconc} will use an alternative formulation of the $\alpha
$-HoC model. A~key observation is that if one generates fitnesses
according to $\alpha$-HoC but then removes interior vertices
independently with some probability $\delta$,  then this results in
a model equivalent to $\alpha'$-HoC for some $\alpha' > \alpha$. The
intuition is that if $\alpha$ is far below the threshold $\frac{\ln
n}{n}$, then not only is there an accessible path with probability
$1-o_n(1)$, but even if we remove a sufficient amount of vertices so
that most paths become forbidden, we will still be below the threshold
and so will still have accessible paths with probability $1-o_n(1)$.
This intuitively requires the original number of accessible paths to be
large. Interestingly, this argument only requires the first equation in
Theorem \ref{thm:aconstr} even though the corollary itself is a
stronger form of that statement.

This idea is formalized in the following lemmas:
%
\begin{lemma}\label{lemma:alphaequivalent2}
Let $\alpha, \delta\in[0, 1]$. Consider the fitness model that first
assigns fitnesses as in $\alpha$-HoC, but then independently removes
each vertex in $\Q_n \setminus\{\bmm{v}^0, \bmm{v}^1\}$ with
probability $\delta$. Then the probability of accessible paths using
only the remaining vertices is $P(n, 1-(1-\alpha)(1-\delta))$.
\end{lemma}
\begin{pf}
Let $\alpha' = 1-(1-\alpha)(1-\delta)$. We compare the model described
above with $\alpha'$-HoC.

Let us make the slight modification to $\alpha'$-HoC and the above
model that we additionally consider any vertex removed if it is less
fit than $\bmm{v}^0$. As no such node can be part of an accessible
path, this will not change accessibility in either model.
We see that these formulations are equivalent up to a translation and
scaling, so they will have the same distribution of accessible paths.
\end{pf}

\begin{lemma}\label{lemma:janson}
Let $\Omega$ be a finite universal set, and let $R$ be a random subset
of $\Omega$ given by $\mathbb{P}(r \in R) = p_r$,
these events being mutually independent over $r\in\Omega$. Let $\{
A_i\}
_{i\in I}$ be subsets of $\Omega$, I a finite
index set. Let $B_i$ be the event $A_i \subseteq R$. Then
%
\begin{equation}
\prod_{i\in I}\mathbb{P}(\bar{B}_i) \leq
\mathbb{P} \biggl( \bigwedge_{i\in I}
\bar{B}_i \biggr).
\end{equation}
\end{lemma}
This inequality is commonly used as a lower bound in Janson's
inequality. See, for instance, Theorem 8.1.1 in \cite{AS}.


\begin{pf*}{Proof of Corollary \ref{cor:aHoCconc}}
The upper bound is simply Markov's inequality. We now turn to the lower
bound. To simplify calculations we may, without loss of generality,
assume that $w_n=o(n\epsilon_n)$ and that $1 \leq w_n \leq
e^{n\epsilon
_n}$ for all $n$.

Let $\delta_n = \epsilon_n-\frac{\ln w_n}{n}$ and let $Y$ denote the
number of intact accessible paths using the same fitness function as
for $X$ but after removing each node except $\bmm{v}^0$ and $\bmm{v}^1$
independently with probability $\delta_n$. By assumption, we know that
$0 \leq\delta_n \leq\epsilon_n \leq\frac{\ln n}{n}$, so $\delta_n$
is always a valid probability.

Using Lemma \ref{lemma:alphaequivalent2} we see that $\mathbb{P}(Y>0) =
P(n, \alpha'_n)$ where $\alpha'_n = 1-(1-\alpha)(1-\delta_n) =
\frac
{\ln n}{n}-\frac{o(1) + \ln w_n}{n}$. As $o(1)+\ln w_n \rightarrow
\infty$ as $n\rightarrow\infty$ it follows from Theorem \ref
{thm:aconstr} that $\lim_{n\rightarrow\infty} \mathbb{P}(Y=0) = 0$.

Condition on the set of accessible paths before removing vertices. Let
$I$ be the set of accessible paths, $R$ the random set of nonremoved
vertices and $B_i$ the event that path $i\in I$ only consist of
nonremoved vertices. Then we are in the setting of Lemma \ref
{lemma:janson}. As the probability that each accessible path remains
intact is $(1-\delta_n)^{n-1}$, averaging conditioned on $X$ we get the
inequality
%
\begin{equation}
\mathbb{P}(Y=0 \mid X) \geq \bigl( 1- (1-\delta_n)^{n-1}
\bigr)^X.
\end{equation}
But since $\lim_{n\rightarrow\infty} \mathbb{P}(Y=0) = 0$ and
$ (
1- (1-\delta_n)^{n-1}  )^X=e^{- (1+o(1) )e^{-n\delta
_n}X}$ it follows that $e^{-n\delta_n}X$ must tend to infinity in
probability. To complete the proof we note that $e^{-n\delta_n}X=\frac
{X}{e^{n\epsilon_n}/w_n} \sim\frac{X}{\mathbb{E}[X]/w_n}$.
\end{pf*}
%

\begin{remark}\label{re2.21}
Note that Proposition \ref{prop:EX2} implies that $\operatorname{Var}(X)\sim
3\mathbb{E}[X]^2$ for $\alpha$ in this regime, so no significant
improvement on Corollary \ref{cor:aHoCconc} can be made by a naive
application of Chebyshev's inequality.
\end{remark}

\section{Results for the RMF model}\label{section:rmf}

Let $n \in\mathbb{N}$, and let $\varepsilon= \varepsilon_n$ be some
strictly positive function. Consider the $n$-dimensional hypercube in
which $\bmm{v}^0$ and $\bmm{v}^1$ are present, and where every other
vertex is present with probability $\varepsilon_n$, independently of
all other vertices. Let $Y = Y_{n,\varepsilon_n}$ denote the number of
accessible paths from $\bmm{v}^0$ to $\bmm{v}^1$, where in this model a
path is accessible if Hamming distance from $\bmm{v}^0$ is strictly
increasing and all vertices along the path are present. The following
proposition may be well known, as it can be interpreted in the context
of site percolation on the directed hypercube. However, we were not
able to locate a suitable reference.

\begin{proposition}\label{vertexremovalprop}
\textup{(i)} $\mathbb{E}[Y] = n! \cdot\varepsilon_n^{n-1}$.

\textup{(ii)} Let $n \rightarrow\infty$, and suppose that $n \varepsilon_n
\rightarrow\infty$. Then $\operatorname{Var}(Y) = o(\mathbb{E}[Y]^2)$, and hence
%
\begin{equation}
Y \sim\mathbb{E}[Y] \sim\frac{\sqrt{2\pi n}}{\varepsilon_n} \biggl( \frac{n\varepsilon_n}{e}
\biggr)^n.
\end{equation}
\end{proposition}

\begin{pf}
There are $n!$ possible paths in the $n$-hypercube. Each path contains
$n-1$ interior vertices, each of which is present with probability
$\varepsilon_n$. This proves~(i). Set $\mu= \mu_n:= n! \varepsilon
_{n}^{n-1}$. Now suppose $n \varepsilon_n \rightarrow\infty$. Let
$Y_i$ be the indicator of the event that the $i$th increasing path is
accessible, where the paths have been ordered in any way. Fix any path
$i_0$. Then, by a standard second moment estimate (see Section~\ref{section:proofaconstr}),
%
\begin{equation}
\operatorname{Var}(Y) \leq\mu+ n! \cdot\sum_{j \sim i_0} \mathbb
{E}(Y_{i_0} Y_j),
\end{equation}
where the sum is taken over all paths $j$ which intersect the path
$i_0$ in at least one interior vertex. Let $k$ be the number of
intersection points. This leaves $T(n, k+1)$ possibilities for the path
$j$. The paths $i_0$ and $j$ contain a total of $2n-2-k$ different
interior vertices; hence the probability of both being present is
$\varepsilon_{n}^{2n-2-k}$. Hence
%
\begin{equation}
\operatorname{Var}(Y) \leq\mu+ n! \cdot\sum_{k=2}^{n}
T(n, k) \varepsilon _{n}^{2n-1-k} \leq\mu+ \mu^2
\cdot\sum_{k=2}^{n} \frac{T(n, k)}{n! \varepsilon
_{n}^{k-1}}.
\end{equation}
Hence since $\mu\rightarrow\infty$ when $n\varepsilon_n \rightarrow
\infty$, it suffices to show that
%
\begin{equation}
\label{eq:thiswas3.4} \sum_{k=2}^{n}
\frac{T(n, k)}{n! \varepsilon_{n}^{k-1}} = o(1).
\end{equation}
We now follow the same strategy as in Section~\ref{section:proofaconstr}, but the analysis
here is much simpler. Let
$\delta\in(0, 1)$. We divide the sum in \eqref{eq:thiswas3.4} into
two parts,~one~for $k \leq(1-\delta) n$ and the other~for $k >
(1-\delta) n$. From Proposition \ref{prop:boundTnksmallk}
and~Lebes\-gue's dominated convergence theorem, it follows easily that, for
any $\delta> 0$, the sum over terms $k \leq(1-\delta)n$ is bounded by
$(1+o_{n}(1)) \sum_{k=2}^{\infty} \frac{k}{(n\varepsilon_n)^{k-1}} =
O(\frac{1}{n\varepsilon_n}) = o(1)$, provided $n\varepsilon_n
\rightarrow\infty$. Similarly, from Proposition \ref
{prop:boundTnklargek} it follows that the sum over terms $k > (1-\delta
)n$ is bounded by
%
\begin{equation}
\label{eq:thiswas3.5} \frac{c}{\mu} \sum_{l=0}^{\delta n}
(l+1) \biggl( \frac{1+2\delta}{5} \cdot n\varepsilon_n
\biggr)^l,
\end{equation}
where $c$ is an absolute constant. Since $n\varepsilon_n \rightarrow
\infty$, the sum in \eqref{eq:thiswas3.5} is bounded by $1+o(1)$ times
the last term, and hence is $O((n\varepsilon_n)^{\delta n})$, which is
in turn $o(\mu)$. This proves \eqref{eq:thiswas3.4} and completes the
proof of the proposition.
\end{pf}

We now turn to the RMF model and prove Theorem \ref{thm:rmf}.

We shall abuse notation and also use $\eta$ to denote the p.d.f. of the
probability distribution under consideration. So suppose $\eta$ has
connected support and is continuous there. Let $\delta> 0$ be given.
Then there exists a bounded, closed interval $I = I_{\delta} \subseteq
\operatorname{Supp}(\eta)$ such that $\int_{I_{\delta}} \eta(x) \,dx >
1-\delta
$. The quantity $c_{\eta,\delta} = \min_{x \in I_{\delta}} \eta(x)$
exists, is nonzero and, obviously, depends only on $\eta$ and $\delta
$. Now let $n \in\mathbb{N}$ and $\theta= \theta_n > 0$ be given.
Without loss of generality, we may assume that the interval $I_{\delta
}$ has length $l(I_{\delta}) > \theta_n / 2$ (in fact any multiple
$c\theta_n$, where $0 < c < 1$, would do in the argument that follows).
By definition of $I_{\delta}$, with probability at least $(1-\delta)^2$
each of $\eta(\bmm{v}^0)$ and $\eta(\bmm{v}^1)$ lie in $I_{\delta}$.
Let $X_{\delta,n,\theta_n}$ be the number of accessible paths in the
$n$-hypercube, where fitnesses are assigned as in \eqref{eq:defrmf},
and conditioning on the fact that both $\eta(\bmm{v}^0)$ and $\eta
(\bmm
{v}^1)$ lie in $I_{\delta}$. We claim that, if $n$ is sufficiently
large, then $X_{\delta,n,\theta_n}$ stochastically dominates the random
variable $Y_{n,\varepsilon_n}$ in Proposition \ref{vertexremovalprop},
where $\varepsilon_n = c_{\eta,\delta} \cdot\frac{\theta_n}{2}$.

To see this, first note that, as long as $l(I_{\delta}) > \theta
_n/2$ then, for any point $x \in I_{\delta}$, there will be an interval
$I_x$ of length at least $\theta_n /2$, which contains $x$ and lies
entirely within $I_{\delta}$. By assumption, any such interval captures
at least $c_{\eta,\delta} \cdot\frac{\theta_n}{2}$ of the distribution
$\eta$. For any adjacent pair $(v, v^{\prime})$ of vertices in the
hypercube such that $d(v^{\prime}, \bmm{v}^0) = d(v, \bmm{v}^0) + 1$,
if $\eta(v^{\prime}) > \eta(v) - \theta_n$, then $v^{\prime}$ is
accessible from $v$. Assuming $\eta(\bmm{v}^0) \in I_{\delta}$, it
follows that we can choose, for each layer $i$ in the hypercube, an
interval $I_i \subseteq I_{\delta}$ of length $\theta_n /2$ such that
any path
%
\begin{equation}\label{eq3.6}
\bmm{v}^0 \rightarrow v_1 \rightarrow v_2
\rightarrow\cdots\rightarrow v_{n-1}
\end{equation}
for which $\eta(v_i) \in I_i$ for all $i = 1, \ldots, n-1$, is
accessible. If $n$ is sufficiently large, we can also ensure that the
interval $I_{n-1}$ contains $\eta(\bmm{v}^1)$, so that any viable path
(\ref{eq3.6}) can definitely be continued to $\bmm{v}^1$. The stochastic
domination of $Y_{n,\varepsilon_n}$ by $X_{\delta,n,\theta_n}$ now
follows. Then one just needs to apply Proposition \ref
{vertexremovalprop} and Theorem~\ref{thm:rmf} follows immediately.

\begin{remark}\label{re3.2}
Suppose $\operatorname{Supp}(\eta)$ is also bounded and that $\theta$ is a constant,
independent of $n$. Let
%
\begin{equation}
C_{\eta,\theta}:= \min_{l(I) = \theta/2, I \subseteq\overline
{\operatorname{Supp}(\eta)}} \int_{I}
\eta(x) \,dx,
\end{equation}
where $I$ denotes a closed interval. Then this minimum exists and is
nonzero. It follows from Proposition \ref{vertexremovalprop} and the
argument above that the number $X = X(n)$ of accessible paths in this
case satisfies
%
\begin{equation}
X \gtrsim n! \cdot C_{\eta,\theta}^{n-1}.
\end{equation}
The point is that $C_{\eta,\theta} \in(0,1]$ is a constant depending
only on $\eta$ and $\theta$.
\end{remark}

\section*{Acknowledgements}

We thank Joachim Krug for making us aware of the problems studied here,
and both he and Stefan Nowak for helpful discussions. We thank both
referees for their very careful reading of the manuscript.

%



\printaddresses

\end{document}